\documentclass[english,a4paper,10pt,fleqn]{article}
\usepackage{theorem}
\usepackage{amssymb,amsmath}
\usepackage{latexsym}
\usepackage{babel}
\usepackage{color}
\usepackage[latin1]{inputenc}
\usepackage{color,graphicx,hyperref,array,fleqn}

\theoremheaderfont{\normalfont\bfseries\upshape}
\theorembodyfont{\normalfont\mdseries\slshape}
\newtheorem{theorem}{Theorem}[section]

\newtheorem{corollary}[theorem]{Corollary}
\newtheorem{lemma}[theorem]{Lemma}
\theorembodyfont{\normalfont\mdseries\upshape}

\newtheorem{definition}[theorem]{Definition}

\numberwithin{equation}{section}

\numberwithin{figure}{section}

\begin{document}

\noindent  \textcolor[rgb]{0.00,0.00,0.60}{\textsf{{\Huge A
theorem of the alternative with an arbitrary number of inequalities and
quadratic programming}}}

\noindent \textbf{M. Ruiz Gal\'an}
\newline \noindent {\small \textsl{University of Granada,
E.T.S. Ingenier\'{\i}a de Edificaci\'on, Department of Applied Mathematics,
c/ Severo Ochoa s/n, 18071 Granada (Spain),
\textcolor[rgb]{0.00,0.00,0.60}{mruizg@ugr.es}}}

\vspace{0.5cm}

\noindent {\small \textbf{Abstract.} In this paper we are concerned with a Gordan-type theorem involving an arbitrary number of inequality functions. We not only state its validity under a weak convexity assumption on the functions, but also show it is an optimal result.
We discuss generalizations of several recent results on nonlinear quadratic optimization, as well as a formula for the Fenchel conjugate of the supremum of a family of functions, in order to illustrate the applicability of that theorem of the alternative.}

\vspace{0.5cm}

\noindent{\small\texttt{\textbf{
\begin{enumerate}
\item Introduction
\item Gordan's theorem for an arbitrary number of inequalities
\item Quadratic optimization 
\item Fenchel conjugate
\end{enumerate}
}}}

\vspace{0.5cm}

\noindent {\small \textbf{2010 Mathematics Subject
Classification:} 90C46, 90C20, 46A22, 26B25.}

\noindent {\small \textbf{Key words:} Theorems of the alternative, quadratic programming, Separation theorem, infsup-convexity.}

\vspace{0.5cm}

\section{Introduction}\label{sect1}

\vspace{-0.5cm}

\textcolor[rgb]{0.89,0.39,0.02}{\textsl{Theorems of the
alternative}} are nothing more than idiosyncratic equivalences
between one statement and the negation of another one --the
alternatives, in most cases in a finite dimensional context,
in which case they are usually reformulations of the Separation
theorem of convex sets in a Euclidean space. However, some of them
predate that fundamental result, as the
\textcolor[rgb]{0.89,0.39,0.02}{\textsl{Gordan theorem}}
\cite{gor}, which dates back to 1873. Notwithstanding its simplicity, Gordan's theorem and some of its extensions are essential tools and active research topics in mathematical
programming (some illustrative examples may be found in
\cite{manbis,jey,dax-sre,gio-gue-thi,jin-kal,jey-lee-li,jey-lee-libis,chu,flo-flo-ver,hu-hua,tan-zha,kur-lee})
and have wide-ranging applications in many different fields (see,
for instance, \cite{gol-evt,li,car-wer,noo-lew-mil,mar}).

One interesting and popular generalization of the Gordan theorem
was given by K. Fan, I. Glicksberg and A.J. Hoffman in \cite{fan-gli-hof}, where an alternative was
established in order to effectively characterize the existence of
a solution for a system with a finite number of convex
inequalities. In \cite{rui2} that result is extended to
inequalities that are convex in a weak sense,
\textcolor[rgb]{0.89,0.39,0.02}{\textsl{infsup-convexity}} (see
Definition \ref{de:infsup} below); in fact, it is shown that
this generalization is optimal, that is, infsup-convexity is the
adequate concept of convexity for dealing with the Gordan theorem.
This is done by means of an equivalent version of the finite
dimensional Separation theorem in the form of minimax inequality. The
main purpose of the present work is to derive a Gordan's theorem
which, on the one hand, is applicable to an arbitrary number of
inequalities and, on the other hand, like its special finite case
\cite{rui2}, is sharp for the kind of convexity under
consideration, again infsup-convexity.

The paper is organized as follows. In Section \ref{sect2} we
deduce, as a consequence of a
Hahn--Banach-type result, the so-called \textcolor[rgb]{0.89,0.39,0.02}{\textsl{Mazur--Orlicz theorem}} (see \cite[Th\'eor\`{e}me 2.41]{maz-orl}, \cite[Theorem, p.
365]{pta}, \cite[Satz, p. 482]{ko02}, \cite[Theorem 28]{sim2} and
its extensions \cite[Theorem 1.1]{ko0}, \cite[Theorem 2.9]{sim1},
\cite[Theorem 12]{grz-prz-urb} and \cite[Theorem 3.1]{din-mo}),
the aforementioned  version of the Gordan theorem for an
arbitrary number, finite or infinite, of inequalities; and at the
same time we establish an equivalence between its validity and the kind of convexity involved, infsup-convexity. The Gordan theorem,
while interesting in its own right, also allows us to
recover and even extend many results on nonlinear optimization. Such is the case on quadratic programming,  as shown in Section
\ref{sect3} regarding several statements by V. Jeyakumar, G.M. Lee and G.Y. Li in \cite{jey-lee-li} on
the solvability of a wide class of quadratic programs, in the form of Karush--Kuhn--Tucker and Fritz Jonh results. In addition, in Section \ref{sect4} we apply the Gordan theorem to dealing with formulae for the conjugate of the supremum of a
family of (possibly infinitely many) functions satisfying no
topological condition.

\section{Gordan's theorem for an arbitrary number of inequalities}\label{sect2}

\vspace{-0.5cm}

This section is devoted to the discussion of the main result in this paper, a version of Gordan's theorem for infinitely many inequality functions. In addition,
we show that such a theorem of the alternative is optimal.

Let us start by recalling that the classical Gordan theorem states that, given $N,m \ge 1$ and $\mathbf{x}_1,\dots, \mathbf{x}_m \in
\mathbb{R}^N$, exactly one of the following alternatives holds:
\begin{enumerate}
\item[{\rm (a1)}] The system
\[
\mathbf{x}_j^T \mathbf{x} < 0  \hbox{ for } j=1,\dots,m
\]
admits a solution $\mathbf{x} \in \mathbb{R}^N$.
\item[{\rm (a2)}] The system
\[
\sum_{j=1}^m t_j \mathbf{x}_j = \mathbf{0}
\]
has a solution $\mathbf{t} \in \Delta_m$,
\end{enumerate}
where $\Delta_m$ is the
\textcolor[rgb]{0.89,0.39,0.02}{\textsl{probability simplex}} in
$\mathbb{R}^m$, that is,  adopting the convention that
$\mathbf{t}$ is the vector $(t_1,\dots,t_m)\in \mathbb{R}^m$,
\[
\Delta_m:= \left\{ \mathbf{t} \in \mathbb{R}^m : \  t_1,\dots,t_m \ge 0
\hbox{ and } \sum_{j=1}^m t_j=1 \right\}.
\]
Worth mentioning is the following generalization of Gordan's theorem, due to K. Fan, I. Glicksberg and A.J. Hoffman (\cite[Theorem 1]{fan-gli-hof}): if $E$ is a real vector space, $C$ is a nonempty and convex subset of $E$, $m\ge 1$ and $f_1,\dots,f_m : C \longrightarrow
\mathbb{R}$ are convex functions, then either
\begin{enumerate}
\item[{\rm (a1)}]
\[
\hbox{there exists } x \in C : \ \max_{j=1,\dots,m} f_j(x) < 0
\]
\end{enumerate}
or
\begin{enumerate}
\item[{\rm (a2)}]
\[
\hbox{for some }\mathbf{t} \in \Delta_m , \   \inf_{x\in C} \sum_{j=1}^m t_j f_j(x) \ge 0,
\]
\end{enumerate}
but never both.

In \cite[Theorem 2.3]{rui2} we have extended this result for
functions satisfying a not very restrictive concept of convexity,
the so-called \textcolor[rgb]{0.89,0.39,0.02}{\textsl{infsup-convexity}} (see Definition \ref{de:infsup}),
and we have shown it is sharp (\cite[Theorem 2.6]{rui2}), in the
sense that the validity of this result implies the
infsup-convexity of the involved functions. Now we follow suit for
an arbitrary number of functions, applied in Section
\ref{sect3} to the study of certain quadratic programs.

But beforehand, let us recall that
given a nonempty set $\Lambda$, $\ell^\infty(\Lambda)$ denotes the
real Banach space of those real-valued functions defined on
$\Lambda$ which are bounded, endowed with its usual addition and
scalar multiplication, as well as with its usual
\textcolor[rgb]{0.89,0.39,0.02}{\textsl{sup norm}}:
\[
\| \Phi \| := \sup_{\lambda \in \Lambda} | \Phi (\lambda )|, \qquad (\Phi  \in \ell^\infty (\Lambda )).
\]
Let us also denote by $\ell^\infty (\Lambda)_+$ the \textcolor[rgb]{0.89,0.39,0.02}{\textsl{nonnegative cone}} of $\ell^\infty (\Lambda )$
\[
\ell^\infty (\Lambda)_+:=\left\{ \Phi \in \ell^\infty (\Lambda) :
\ 0 \le \inf_{\lambda \in \Lambda} \Phi (\lambda ) \right\}.
\]
In what follows, we identify a function $\Phi : \Lambda \longrightarrow \mathbb{R}$ with its values $\{ \Phi(\lambda) \}_{\lambda \in \Lambda}$, and write $\| \cdot \|$ for the sup norm and $\|
\cdot \|_*$ for its dual norm. In addition, $\mathbf{1}$ stands for the function in
$\ell^\infty (\Lambda )$ defined at each $\lambda \in \Lambda$ as
\[
\mathbf{1}(\lambda):=1.
\]

We begin by introducing a result which is elementary but useful for our purposes. We include its proof so that our approach will be as self-contained as possible.
Let us first recall that if $\Lambda$ is a nonempty set, a function $F : \ell^\infty (\Lambda) \longrightarrow \mathbb{R}$ is said to be \textcolor[rgb]{0.89,0.39,0.02}{\textsl{positive}} provided that
\[
\Phi \in \ell^\infty (\Lambda)_+ \ \Rightarrow \ F (\Phi ) \ge 0.
\]

\bigskip

\begin{lemma}\label{le:functionals}
Let $\Lambda$ be a nonempty set and $L : \ell^\infty (\Lambda) \longrightarrow \mathbb{R}$ be a linear functional. Then the following assertions are equivalent:
\begin{enumerate}
\item[{\rm (i)}]  $\Phi \in \ell^\infty(\Lambda) \ \Rightarrow \ L(\Phi ) \le \displaystyle \sup_{\lambda \in \Lambda} \Phi (\lambda).$

\item[{\rm (ii)}] $L \hbox{ is positive and } L(\mathbf{1})=1$.
\end{enumerate}
Moreover, if some of these equivalent statements hold, then $L$ is continuous and $\|L\|_*=1$.
\end{lemma}

\noindent \textsc{Proof.} \fbox{(i) $\Rightarrow$ (ii)}
$L$ is positive, since for all $\Phi
\in \ell^\infty (\Lambda)_+$, 
\[
\begin{array}{rl}
  L(-\Phi ) & \le \displaystyle \sup_{\lambda \in \Lambda} -\Phi (\lambda)   \\
   & \le 0,
\end{array}
\]
hence $L(\Phi ) \ge 0$.
On the other hand,
\[
\left. \begin{array}{rll}
         L(\mathbf{1}) & \le \displaystyle \sup_{\lambda \in \Lambda} 1 & =1 \\
         -L(\mathbf{1}) & \le \displaystyle \sup_{\lambda \in \Lambda} -1  & =-1
       \end{array}
 \right\} \ \Rightarrow \ L(\mathbf{1})=1.
\]

\noindent \fbox{(ii) $\Rightarrow$ (i)} Given $\Phi \in
\ell^\infty (\Lambda)$, $\left( \sup_{\lambda \in \Lambda} \Phi (\lambda)
\right)\mathbf{1}-\Phi \in \ell^\infty (\Lambda)_+$, so, taking into account that $L$ is positive, linear and
$L(\mathbf{1})=1$, we arrive at
\[
L(\Phi ) \le \sup_{\lambda \in \Lambda} \Phi (\lambda).
\]
The arbitrariness of $\Phi \in \ell^\infty(\Lambda)$ yields (i).

Finally, for each $\Phi \in \ell^\infty (\Lambda)$ there holds
\[
\begin{array}{rl}
  L(\Phi ) & \le \displaystyle  \sup_{\lambda \in \Lambda} \Phi (\lambda) \\
   & \le \| \Phi \|,
\end{array}
\]
hence, the linear functional $L$ is continuous and $\|L\|_* \le 1$, and since $L(\mathbf{1})=1$, then $\|L\|_*=1$.
\hfill$\textcolor[rgb]{0.00,0.00,0.60}{\Box}$

When $\Lambda$ is finite, $\ell^\infty (\Lambda)= \mathbb{R}^m$
for some $m \ge 1$. 
Then $L: \mathbb{R}^m \longrightarrow \mathbb{R}$ linear, positive with $L(\mathbf{1})=1$ means $L=\mathbf{t}$ for some $\mathbf{t} \in \Delta_m$. This, together with Lemma \ref{le:functionals}, motivates the introduction of the following notation: given a nonempty set $\Lambda$, we write
\[
\Delta_\Lambda := \{ L : \ell^\infty (\Lambda) \longrightarrow \mathbb{R}: \ L \hbox{ is positive, linear and } L(\mathbf{1})=1 \}.
\]

We are going to deduce our Gordan-type result from a generalization of the Hahn--Banach theorem, the
\textcolor[rgb]{0.89,0.39,0.02}{\textsl{Mazur--Orlicz
theorem}} (\cite[Th\'eor\`{e}me 2.41]{maz-orl}, \cite[Theorem, p.
365]{pta}, \cite[Satz, p. 482]{ko02}, \cite[Theorem 1.1]{sim2}), that asserts that if $E$ is a real vector space, $C$ is a nonempty and convex subset
of $E$, and $S: E \longrightarrow \mathbb{R}$ is a sublinear (subadditive and positively homogeneous)
functional, then, there exists a linear functional $L: E
\longrightarrow \mathbb{R}$ such that
\[
x \in E \ \Rightarrow \ L(x) \le S(x)
\]
and
\[
\inf_{x \in C} L(x) = \inf_{x \in C} S(x).
\]
We also need to recall the aforementioned concept of convexity, infsup-convexity, that arises in minimax theory (see \cite[p. 653]{kas-kol} and \cite[Definition
2.11]{ste}), where it has turned out to be suitable (\cite[Theorem
2.20]{rui3} and \cite[Corollary 3.12]{rui2}), as in the finite version of Gordan's theorem (\cite[Theorem 2.6]{rui2}).

\bigskip

\begin{definition}\label{de:infsup}
Given $\Lambda$ and $X$ nonempty sets, a family of real valued functions on $X$, $\left\{ f_\lambda
\right\}_{\lambda \in \Lambda}$, is said to be \textcolor[rgb]{0.89,0.39,0.02}{\textsl{infsup-convex on}} $X$
provided that
\[
\left. \begin{array}{c}
         m\ge 1, \ \mathbf{t} \in \Delta_m \\
         x_1,\dots , x_m \in X
       \end{array}
 \right\}
 \ \Rightarrow \
 \displaystyle \inf_{x\in X}\sup_{\lambda\in \Lambda} f_\lambda (x) \le \sup_{\lambda\in \Lambda}
\sum_{j=1}^m t_j f_\lambda (x_j).
\]
\end{definition}

\bigskip

Likewise we can define
\textcolor[rgb]{0.89,0.39,0.02}{\textsl{supinf-concavity}} of $f$
on $X$, although this notion is not need in this paper
(see \cite[Definition 2.1]{rui2}).

For instance, if $X$ is a nonempty and convex subset of a real vector space, $\Lambda$ is a nonempty set and, for each $\lambda \in \Lambda$, $f_\lambda : X \longrightarrow
\mathbb{R}$ is a convex function, then the family $\left\{
f_\lambda \right\}_{\lambda \in \Lambda}$ is infsup-convex on $X$.
We can also check easily that if $X$ and $\Lambda$
are nonempty sets and $\left\{ f_\lambda
\right\}_{\lambda \in \Lambda}$ is a family of real-valued
functions defined on $X$ such that
\[
x \in X \ \Rightarrow \ \left\{ f_\lambda (x) \right\}_{\lambda
\in \Lambda} \in \ell^\infty (\Lambda)
\]
and the joint range set
\[
\left\{  \left\{ f_\lambda (x) \right\}_{\lambda \in \Lambda} : \
x \in X \right\} \subset \ell^\infty (\Lambda)
\]
is convex, then the family is infsup-convex on $X$. Some well-known results guarantee the convexity of the joint range set. Thus, if for $N \ge 1$, $\mathbb{S}^N$ denotes the set of the \textcolor[rgb]{0.89,0.39,0.02}{\textsl{$N \times N$ symmetric real matrices}}, then the Dines theorem (\cite[Theorem 1]{din}) asserts that for any $\mathbf{A}_1,\mathbf{A}_2 \in \mathbb{S}^N$ the joint range set of the corresponding quadratic forms,
\[
\left\{ \left( \frac{1}{2} \mathbf{x}^T \mathbf{A}_1 \mathbf{x},\frac{1}{2} \mathbf{x}^T \mathbf{A}_2 \mathbf{x} \right) : \ \mathbf{x} \in \mathbb{R}^N \right\},
\]
is a convex subset of $\mathbb{R}^2$. In this respect, it is worth mentioning the Brickman theorem \cite[Theorem 2.1]{bri}, which asserts that the subset of that convex set
\[
\left\{ \left( \frac{1}{2} \mathbf{x}^T \mathbf{A}_1 \mathbf{x},\frac{1}{2} \mathbf{x}^T \mathbf{A}_2 \mathbf{x} \right) : \ \mathbf{x} \in \mathbb{R}^N, \ \| \mathbf{x} \|=1 \right\},
\]
is also convex as soon as $N \ge 3$. Although these results fail for three (or more) quadratic forms, there
are sufficient conditions for the convexity of the joint range set
of three such forms: see, for instance, the Polyak theorem
\cite[Theorem 2.1]{pol} and its generalization \cite[Theorem
4.1]{jey-lee-li}.

We now come to our main result, a
theorem of the alternative of the Gordan-type generalizing the
finite case \cite[Theorem 3.2]{rui2}.

\bigskip

\begin{theorem}\label{th:gordan}
Suppose that $\Lambda$ and $X$ are nonempty sets and that $\left\{ f_\lambda  \right\}_{\lambda \in \Lambda}$ is a family of real valued functions on $X$ that is infsup-convex on $X$ and satisfying
\[
x \in X \ \Rightarrow \ \left\{ f_\lambda (x)
\right\}_{\lambda \in \Lambda} \in \ell^\infty (\Lambda).
\]
Then, one, and only one, of the following statements holds:
\begin{enumerate}
\item[{\rm (a1)}] There exists $x \in X$ such that
\[
\sup_{\lambda \in \Lambda}  f_\lambda(x) <0.
\]
\item[{\rm (a2)}] There exists $L \in \Delta_\Lambda$ such that
\[
\inf_{x \in X} L\left( \left\{ f_\lambda (x) \right\}_{\lambda \in \Lambda} \right) \ge 0.
\]
\end{enumerate}
\end{theorem}

\noindent \textsc{Proof.} Consider the real linear vector space $\ell^\infty (\Lambda)$, the
nonempty and convex subset of $\ell^\infty (\Lambda)$
\[
C:=\textrm{conv} \left\{ \{f_\lambda (x) \}_{\lambda \in \Lambda}: \
x \in X\right\},
\]
(``conv" denotes ``convex hull"), and the sublinear functional $S: \ell^\infty (\Lambda) \longrightarrow \mathbb{R}$ defined at each $\Phi \in \ell^\infty (\Lambda)$ as $\displaystyle S(\Phi):=\sup_{\lambda \in \Lambda} \Phi(\lambda)$. Then, the Mazur--Orlicz theorem and Lemma \ref{le:functionals} imply the existence of $L \in
\Delta_\Lambda$ such that
\[
\inf_{\Phi \in C} L(\Phi )= \inf_{\Phi \in C} S(\Phi).
\]
In particular, we have the alternative:
\[
\inf_{\Phi \in C} L(\Phi ) \ge 0
\]
or (exclusive)
\[
\inf_{\Phi \in C} S(\Phi)<0.
\]
This completes the proof, since
\[
\inf_{\Phi \in C} L(\Phi ) = \inf_{x \in X} L\left( \{f_\lambda (x)\}_{\lambda \in \Lambda}\right),
\]
while, according to the infsup-convexity of the family $\left\{ f_\lambda \right\}_{\lambda \in \Lambda}$ on $X$, 
\[
\begin{array}{rl}
\displaystyle \inf_{\Phi \in C} \sup_{\lambda \in \Lambda}
\Phi (\lambda) & = \displaystyle \inf_{\substack{n\ge 1, \ \mathbf{s} \in \Delta_n \\
x_1,\dots,x_n\in X }} \sup_{\lambda \in \Lambda} \sum_{i=1}^n s_i f_\lambda (x_i ) \\
&  \\
&  = \displaystyle \inf_{x \in X} \sup_{\lambda \in \Lambda} f_\lambda (x).
\end{array}
\]
\hfill$\textcolor[rgb]{0.00,0.00,0.60}{\Box}$

\bigskip

In some theorems of the alternative for quadratic inequalities, the convexity of the joint range set is
proven and then the Separation theorem is suitably applied, as in the celebrated
\textcolor[rgb]{0.89,0.39,0.02}{\textsl{S-Lemma}}, established by
V.A. Yakubovi\v{c} from the Dines theorem \cite{yak}, or in its
generalization \cite[Corollary 3.7]{jey-lee-li}. Instead,
sometimes one looks for a  convex set containing the joint range
set and applies the Separation theorem: see for instance
\cite[Theorem 2.4]{hu-hua}. We will show in the next section that
the use of the theorem of Gordan, Theorem \ref{th:gordan}, avoids
this kind of technique and leads to a more direct approach for
the study of quadratic programs, which in particular allows us to
deduce quite general results. As an advance of it, let us consider
the following corollary which, in the finite case and when
$X=\mathbb{R}^N$, coincides with Yuan's alternative theorem
\cite[Lemma 2.3]{yua}:

\bigskip

\begin{corollary}
Let $N \ge 1$, $\mathbf{A}_1,\mathbf{A}_2 \in \mathbb{S}^N$ let
$X$ be $\mathbb{R}^N$ or $ \{ \mathbf{x} \in \mathbb{R}^N : \ \|
\mathbf{x} \| =1 \}$. Then, either
\begin{enumerate}
\item[{\rm (a1)}] there exists $\mathbf{x} \in X$ such that
\[
\max_{i=1,2} \left\{ \frac{1}{2} \mathbf{x}^T \mathbf{A_i} \mathbf{x} \right\} <0
\]
\end{enumerate}
or
\begin{enumerate}
\item[{\rm (a2)}] there exists $0 \le t \le 1$ in such a way that the matrix $t\mathbf{A}_1+(1-t)\mathbf{A}_2$ is semi-definite positive,
\end{enumerate}
but never both.
\end{corollary}

\noindent \textsc{Proof.} Since the joint range set is convex, according to Dines' theorem \cite[Theorem 1]{din} and Brickman's theorem \cite[Theorem 2.1]{bri}, then the family made up of the functions $f_1,f_2 : X \longrightarrow \mathbb{R}$ defined at each $\mathbf{x} \in X$ as
\[
f_i(\mathbf{x}):= \frac{1}{2}\mathbf{x}^T \mathbf{A}_i \mathbf{x}, \qquad (i=1,2)
\]
is infsup-convex on $X$, and therefore, it suffices to apply the Gordan theorem, Theorem \ref{th:gordan}.
\hfill$\textcolor[rgb]{0.00,0.00,0.60}{\Box}$

\bigskip

The next direct consequence of Theorem \ref{th:gordan} is the
extension of the original Gordan theorem to an arbitrary number of vectors:

\bigskip

\begin{corollary} Assume that $E$ is a real normed space, $\Lambda$ is a nonempty set and that for all $\lambda \in \Lambda$, $x_\lambda^* \in E^*$, in such a way that the subset of $E^*$ $\{ x_\lambda ^* : \ \lambda \in \Lambda \}$ is bounded. Then the problem
\[
\hbox{find } x \in E  \hbox{ such that } \sup_{\lambda \in \Lambda} x_\lambda^* (x) < 0
\]
is solvable if, and only if, its dual problem
\[
\hbox{find } L \in \Delta_\Lambda \hbox{ such that } L\left( \{ x_\lambda^* (\cdot )\}_{\lambda \in \Lambda}\right)=0
\]
has no solution.
\end{corollary}

\noindent \textsc{Proof.}  Apply Gordan's theorem, Theorem
\ref{th:gordan}, with the nonempty set
\[
X:=E
\]
and the family of functions
\[
f_\lambda (x):= x_\lambda^*(x), \qquad (\lambda \in \Lambda, \ x \in X),
\]
which is clearly infsup-convex on $E$ and in addition satisfies
\[
x \in E \ \Rightarrow \ \{ x_\lambda^*(x)\}_{\lambda \in \Lambda} \in \ell^\infty (\Lambda).
\]
Hence, that theorem yields that, either
\[
\hbox{there exists } x \in E : \ \sup_{\lambda \in \Lambda} x_\lambda^* (x) < 0
\]
or (exclusive)
\[
\hbox{there exists } L \in \Delta_\Lambda : \ x \in E \ \Rightarrow \ \inf_{x \in E} L\left( \{ x_\lambda^* (x)\}_{\lambda \in \Lambda}\right) \ge 0.
\]
But this second condition is exactly
\[
\inf_{x \in E} h(x) \ge 0,
\]
with $h : E \longrightarrow \mathbb{R}$ being the linear
functional given for any $x \in E$ as
\[
h(x):=L\left( \{ x_\lambda^* (x)\}_{\lambda \in \Lambda}\right),
\]
and therefore ($h$ is linear)
\[
h=0
\]
or, in other words,
\[
x \in E \ \Rightarrow \ \inf_{x \in E} L\left( \{ x_\lambda^* (x)\}_{\lambda \in \Lambda}\right)=0,
\]
i.e.,
\[
L\left( \{ x_\lambda^* (\cdot )\}_{\lambda \in \Lambda}\right)=0.
\]
\hfill$\textcolor[rgb]{0.00,0.00,0.60}{\Box}$

\bigskip

Notice that, in view of the Separation theorem, the first
condition means that $0$ does not belong to the weak-$^*$ closure
of $\textrm{conv}\{ x_\lambda^*: \ \lambda \in \Lambda \}$.

Now we prove that our theorem of the alternative, Theorem
\ref{th:gordan}, is optimal. More specifically, it is obvious that
the role of the scalar $0$ in (\cite[Theorem 1]{fan-gli-hof} and)
Theorem \ref{th:gordan} is irrelevant, or in other words, the
following version of Theorem \ref{th:gordan} remains true: if
$\alpha \in \mathbb{R}$, $\Lambda$ and $X$ are nonempty sets and
$\left\{ f_\lambda \right\}_{\lambda \in \Lambda}$ is a family
that is infsup-convex on $X$ and satisfies
\[
x \in X \ \Rightarrow \ \left\{ f_\lambda (x) \right\}_{\lambda \in \Lambda} \in \ell^\infty (\Lambda),
\]
then, one, and only one, of the following statements holds:
\begin{enumerate}
\item[{\rm (a1)}] There exists $x \in X$ such that
\[
\sup_{\lambda \in \Lambda}  f_\lambda (x) < \alpha.
\]
\item[{\rm (a2)}] There exists $L \in \Delta_\Lambda$ with
\[
\inf_{x \in X} L\left( \left\{ f_\lambda (x) \right\}_{\lambda \in \Lambda} \right) \ge \alpha
\]
\end{enumerate}
(use that $L(\mathbf{1})=1$ and the elementary fact that $\left\{ f_\lambda \right\}_{\lambda \in \Lambda}$ is infsup-convex on $X$ if, and only if, $\left\{ f_\lambda - \alpha \right\}_{\lambda \in \Lambda}$ is also infsup-convex on $X$).

One of the two hypotheses is clearly necessary in order that the
alternative makes sense (for all $x \in X$, $\left\{ f_\lambda (x)
\right\}_{\lambda \in \Lambda} \ell^\infty(\Lambda)$). The other
one, infsup-convexity on $X$, is also necessary for the validity
of the alternative, as we show in Theorem
\ref{th:characterization} below. Beforehand, an easy remark
extending a well-known result in the finite case (see for instance
\cite[Lemma 2.6]{rui3}).

\bigskip

\begin{lemma}\label{le:elemental}
If $\Lambda$ is a nonempty set, $L \in \Delta_\Lambda$, $m \ge 1$,
$\mathbf{t}\in \Delta_m$ and $\{ \alpha_\lambda^{(1)} \}_{\lambda
\in \Lambda},\dots, \{ \alpha_\lambda^{(m)} \}_{\lambda \in
\Lambda} \in \ell^\infty (\Lambda)$, then
\[
\min_{j=1,\dots,m} L\left( \{ \alpha_\lambda^{(j)} \}_{\lambda \in \Lambda} \right)  \le
\sup_{\lambda \in \Lambda}  \sum_{j=1}^m t_j \alpha_\lambda^{(j)}.
\]
\end{lemma}

\noindent \textsc{Proof.} It is a direct consequence of the
following chain of inequalities:
\[
\begin{array}{rl}
  \displaystyle \min_{j=1,\dots,m} L\left( \{ \alpha_\lambda^{(j)} \}_{\lambda \in \Lambda} \right)  & \displaystyle  \le \sum_{j=1}^m t_j L \left( \{ \alpha_\lambda^{(j)} \}_{\lambda \in \Lambda} \right) \qquad (t \in \Delta_m)  \\
   & \displaystyle = L \left( \sum_{j=1}^m t_j \{ \alpha_\lambda^{(j)} \}_{\lambda \in \Lambda} \right)   \\
   & \displaystyle \le \sup_{\lambda \in \Lambda } \sum_{j=1}^m t_j  \alpha_\lambda^{(j)}  \qquad \hbox{(Lemma \ref{le:functionals})}.
\end{array}
\]
\hfill$\textcolor[rgb]{0.00,0.00,0.60}{\Box}$

\bigskip

Now we are in a position to obtain the announced optimality of the
version of the Gordan theorem --for (finitely or) infinitely many
functions given in Theorem \ref{th:gordan}-- that extends the finite
case stated in \cite[Theorem 2.6]{rui2}:

\bigskip

\begin{theorem}\label{th:characterization}
Let $\Lambda$ and $X$ be nonempty sets and $\left\{ f_\lambda \right\}_{\lambda \in \Lambda}$ be a family of real valued functions on $X$ such that
\[
x \in X \ \Rightarrow \ \left\{ f_\lambda (x) \right\}_{\lambda \in \Lambda} \in \ell^\infty (\Lambda).
\]
Then
\[
\left\{ f_\lambda \right\}_{\lambda \in \Lambda} \hbox{ is infsup-convex on } X
\]
if, and only if, for all $\alpha \in \mathbb{R}$, exactly one of
the following conditions holds:
\begin{enumerate}
\item[{\rm (a1)}] There exists $x \in X$ such that
\[
\sup_{\lambda \in \Lambda}  f_\lambda (x) < \alpha.
\]
\item[{\rm (a2)}] There exists $L \in \Delta_\Lambda$ such that
\[
\inf_{x \in X} L\left( \left\{ f_\lambda (x) \right\}_{\lambda \in \Lambda} \right) \ge \alpha.
\]
\end{enumerate}
\end{theorem}

\noindent \textsc{Proof.} According to Theorem \ref{th:gordan} and
the above discussion, what remains is only to prove the sufficiency.
So, let $m\ge 1$, $\mathbf{t} \in \Delta_m$ and $x_1,\dots,x_m \in
X$. We assume without any loss of generality that
\[
\alpha := \inf_{x \in X} \sup_{\lambda \in \Lambda} f_\lambda (x) >
-\infty .
\]
The alternative (a1) fails, so there exists $L \in \Delta_\Lambda$
with
\begin{equation}\label{eq:ii}
    \inf_{x \in X} L\left( \left\{ f_\lambda (x) \right\}_{\lambda \in \Lambda}\right) \ge \alpha.
\end{equation}
Therefore
\[
\begin{array}{rl}
  \displaystyle  \inf_{x \in X} \sup_{\lambda \in \Lambda} f_\lambda (x) & = \alpha \\
   & \displaystyle \le \inf_{x \in X} L\left( \left\{ f_\lambda (x) \right\}_{\lambda \in \Lambda} \right) \qquad \hbox{(by \eqref{eq:ii})} \\
   &  \le \displaystyle \min_{j=1,\dots, m} L\left( \left\{ f_\lambda (x_j) \right\}_{\lambda \in \Lambda} \right)  \\
   & \displaystyle \le \sup_{\lambda \in \Lambda}  \sum_{j=1}^m t_j f_\lambda (x_j)   \qquad \hbox{(by Lemma \ref{le:elemental})} \\
\end{array}
\]
and so,
\begin{displaymath}
\inf_{x \in X} \sup_{\lambda \in \Lambda} f_\lambda (x) \le
\sup_{x \in X} \sum_{j=1}^m t_j f_\lambda (x_j). \\
\end{displaymath}
As $m\ge 1$, $\mathbf{t} \in \Delta_m$ and $x_1,\dots,x_m \in X$
are arbitrary, then the family $\left\{ f_\lambda
\right\}_{\lambda \in \Lambda}$ is infsup-convex on $X$, as was to be shown.
\hfill$\textcolor[rgb]{0.00,0.00,0.60}{\Box}$

\bigskip

Even in the finite case, this result provides us with examples of
infsup-convex families of quadratic functions. Thus, thanks to the
nonhomogeneous version of the Yuan theorem given in \cite[Theorem
3.6]{jey-lee-li}, if $N \ge 1$, $\mathbf{A}_1,\mathbf{A}_2 \in
\mathbb{S}^N$, $\mathbf{b}_1,\mathbf{b}_2 \in \mathbb{R}^N$ and
$c_1,c_2 \in \mathbb{R}$, then the family of quadratic functions
$\{q_1,q_2\}$ defined at each $\mathbf{x} \in \mathbb{R}^N$ as
\[
q_i(\mathbf{x}):=\frac{1}{2}\mathbf{x}^T \mathbf{A}_i \mathbf{x} + \mathbf{b}_i^T \mathbf{x}+c_i ,  \qquad (i=1,2),
\]
is infsup-convex on any affine subspace of $\mathbb{R}^N$. In a
similar fashion,  we can use Theorem \ref{th:characterization} and
\cite[Theorem 2.4]{hu-hua}, \cite[Theorem 4.2]{jey-lee-li} or
\cite[Corollary 3.5]{jey-lee-li} to obtain examples of
infsup-convex families of quadratic functions. In the next section we
describe a class of (possibly) infinitely many quadratic functions that is infsup-convex on $\mathbb{R}^N$.

\section{Quadratic optimization}\label{sect3}

\vspace{-0.5cm}

We first illustrate the application of the Gordan theorem, Theorem
\ref{th:gordan}, with a theorem of the alternative for quadratic inequalities,
Corollary \ref{co:quadratic}, which generalizes that established
in \cite[Theorem 5.2]{jey-lee-li} for a finite number of quadratic
functions to the infinite case, even for domains more general than
$\mathbb{R}^N$. In particular, we study the solvability of certain
nonlinear, infinite and quadratic programs in terms of suitable
Karush--Kuhn--Tucker and Fritz John conditions, extending those
given in \cite[Corollary 5.3]{jey-lee-li}.

First, let us recall that given $N \ge 1$, a real matrix $\mathbf{A} =\left(
a_{kl}\right)_{k,l=1,\dots,N }$ is said to be a \textcolor[rgb]{0.89,0.39,0.02}{\textsl{Z-matrix}} provided that $A \in \mathbb{S}^N$ and
\[
k<l \ \Rightarrow \ a_{kl} \le 0.
\]
Both this kind of matrix and its generalizations have found many
applications: it suffices to see for instance
\cite{bap-rag,che-xia,lou-qin-kon-xiu,gow-rav}.

We write $\mathbb{R}^N_+:=\{ \mathbf{x} \in \mathbb{R}^N : \ x_1,\dots,x_N \ge 0 \}$.
\bigskip

\begin{lemma}\label{le:quadratic}
Let $N \ge 1$, $\Lambda$ be a nonempty set and suppose that for each
$\lambda \in \Lambda$, $\mathbf{A}_\lambda  \in \mathbb{S}^N$,
$\mathbf{b}_\lambda \in \mathbb{R}^N$, $c_\lambda \in \mathbb{R}$
in such a way that
\[
\left( \begin{array}{cc}
  \mathbf{A}_\lambda & \mathbf{b}_\lambda \\
  \mathbf{b}_\lambda ^T & 2 c_\lambda
\end{array} \right)
\]
is a Z-matrix. If in addition for any $\lambda \in \Lambda$
$q_\lambda : \mathbb{R}^N \longrightarrow \mathbb{R}$ is the
quadratic function defined at each $\mathbf{x} \in \mathbb{R}^N$ as
\[
q_\lambda (\mathbf{x}):= \frac{1}{2}\mathbf{x}^T
\mathbf{A}_\lambda \mathbf{x} + \mathbf{b}_\lambda ^T \mathbf{x} +
c_\lambda,
\]
and $X$ is a subset of $\mathbb{R}^N$ with $\mathbb{R}^N_+ \subset
X$, then the family $\{q_\lambda\}_{\lambda \in \Lambda}$ is
infsup-convex on $X$.
\end{lemma}

\noindent \textsc{Proof.}  Let $m \ge 1$, $\mathbf{t} \in
\Delta_m$ and $\mathbf{x}_1=\left( x_1^{(1)},\dots,x_N^{(1)}
\right) ,\dots,\mathbf{x}_m =\left( x_1^{(m)},\dots,x_N^{(m)}
\right) \in X$. We must prove
\begin{equation}\label{eq:infsup}
\inf_{\mathbf{x}\in X} \sup_{\lambda \in \Lambda} q_\lambda (\mathbf{x}) \le \sup_{\lambda \in \Lambda} \sum_{j=1}^m t_j q_\lambda (\mathbf{x}_j).
\end{equation}
Indeed, we take the element $\mathbf{x}_0 \in X$, which does not depend on
$\lambda \in \Lambda$,
\[
x_k^0:= \sqrt{\sum_{j=1}^m t_j \left( x_k^{(j)} \right)^2}, \qquad (k=1,\dots,N),
\]
--well defined, because $\mathbf{t} \in \Delta_m$, and since
$\mathbb{R}^N_+ \subset X$, it belongs to $X$. Then, given
$\lambda \in \Lambda$, if we write $\mathbf{A}_\lambda=\left(
a_{kl}^{(\lambda)}\right)_{k,l=1,\dots,N }$ and
$\mathbf{b}_\lambda =\left(
b_1^{(\lambda)},\dots,b_N^{(\lambda)}\right)$, we have that
\[
\begin{array}{rl}
q_\lambda (\mathbf{x}_0) & = \displaystyle  \frac{1}{2} \sum_{k=1}^N \left( x_k^{(0)}
        \right)^2  a_{kk}^{(\lambda )} +  \sum_{\substack{k,l=1 \\ k <
l}}^N  x_k^{(0)} x_l^{(0)} a_{kl}^{(\lambda)} + \sum_{k=1}^N
x_k^{(0)} b_k^{(\lambda)} + c_ \lambda \\
	& \displaystyle \le \frac{1}{2} \sum_{k=1}^N \left( x_k^{(0)}
        \right)^2  a_{kk}^{(\lambda )} +  \sum_{\substack{k,l=1 \\ k <
l}}^N  x_k^{(0)} x_l^{(0)} a_{kl}^{(\lambda)} + \sum_{k=1}^N
\left( \sum_{j=1}^m t_j
x_k^{(j)} \right) b_k^{(\lambda)} + c_ \lambda \\
	& \displaystyle \le \frac{1}{2} \sum_{k=1}^N \left( x_k^{(0)} \right)^2  a_{kk}^{(\lambda )} +  \sum_{\substack{k,l=1 \\ k < l}}^N \left( \sum_{j=1}^m t_j x_k^{(j)} x_l^{(j)} \right) a_{kl}^{(\lambda)} + \sum_{k=1}^N \left( \sum_{j=1}^m t_j x_k^{(j)} \right) b_k^{(\lambda)} + c_ \lambda \\
    & \displaystyle = \frac{1}{2} \sum_{k=1}^N \left( \sum_{j=1}^m t_j \left( x_k^{(j)} \right)^2 \right) a_{kk}^{(\lambda )} +  \sum_{\substack{k,l=1 \\ k < l}}^N \left( \sum_{j=1}^m t_j x_k^{(j)} x_l^{(j)} \right) a_{kl}^{(\lambda)} + \sum_{k=1}^N \left( \sum_{j=1}^m t_j x_k^{(j)} \right) b_k^{(\lambda)} + c_ \lambda \\
        & \displaystyle = \frac{1}{2} \sum_{j=1}^m t_j \mathbf{x}_j^T \mathbf{A}_\lambda \mathbf{x}_j + \sum_{j=1}^m t_j \mathbf{b}_\lambda^T \mathbf{x}_j +c_\lambda   \\
	& \displaystyle  = \sum_{j=1}^m t_j q_\lambda (\mathbf{x}_j), 
\end{array}
\]
where in the first inequality we have used the fact that for all $k=1,\dots,N$, $b_k^{(\lambda )} \le 0$ and the
Cauchy--Schwarz inequality for the vectors $\left( \sqrt{t_1}x_k^{(1)}, \dots,
\sqrt{t_m}x_k^{(m)} \right)$ and $\left( \sqrt{t_1}, \dots,
\sqrt{t_m} \right)$; and in the second that 
$\mathbf{A}_\lambda$ is a Z-matrix and again the Cauchy--Schwarz
inequality for $\left( \sqrt{t_1}x_k^{(1)}, \dots,
\sqrt{t_m}x_k^{(m)} \right)$ and $\left( \sqrt{t_1}x_l^{(1)},
\dots, \sqrt{t_m}x_l^{(m)} \right)$. Finally, the arbitrariness of $\lambda \in \Lambda$ and the fact that $\mathbf{x}_0$ belongs to $X$ imply \eqref{eq:infsup}, and we are done.
 
 \hfill$\textcolor[rgb]{0.00,0.00,0.60}{\Box}$

\bigskip

We deduce, as a consequence of Lemma \ref{le:quadratic} and Theorem \ref{th:gordan}, the following alternative for quadratic inequalities associated with Z-matrices, which not only extends the
finite case in \cite[Theorem 5.2]{jey-lee-li} to the infinite one, but also, even in the finite case, allows one to consider sets $X$ more general than $\mathbb{R}^N$:

\bigskip

\begin{corollary}\label{co:quadratic}
Suppose that $N \ge 1$, $X$ is a subset of $\mathbb{R}^N$ such that $\mathbb{R}^N_+ \subset X$, $\Lambda$ is a nonempty set, and that for each $\lambda \in \Lambda$,
$\mathbf{A}_\lambda  \in \mathbb{S}^N$,
$\mathbf{b}_\lambda \in \mathbb{R}^N$,
$c_\lambda \in \mathbb{R}$ are such that
\[
\left( \begin{array}{cc}
  \mathbf{A}_\lambda & \mathbf{b}_\lambda \\
  \mathbf{b}_\lambda ^T & 2 c_\lambda
\end{array} \right)
\]
is a Z-matrix. If, in addition, the quadratic function
$q_\lambda : \mathbb{R}^N \longrightarrow \mathbb{R}$
\[
q_\lambda (\mathbf{x}):= \frac{1}{2}\mathbf{x}^T
\mathbf{A}_\lambda \mathbf{x} + \mathbf{b}_\lambda ^T \mathbf{x} +
c_\lambda, \qquad (\mathbf{x} \in \mathbb{R}^N)
\]
satisfies
\[
\mathbf{x} \in X \ \Rightarrow \ \left\{ q_\lambda (\mathbf{x}) \right\}_{\lambda \in \Lambda} \in \ell^\infty (\Lambda),
\]
then, one, and only one, of the following alternatives holds:
\begin{enumerate}
\item[{\rm (a1)}] There exists $\mathbf{x} \in X$ with
\[
\sup_{\lambda \in \Lambda}  q_\lambda(\mathbf{x}) <0.
\]
\item[{\rm (a2)}] There exists $L \in \Delta_\Lambda$ such that
\[
\inf_{\mathbf{x} \in X} L\left(\left\{ q_\lambda (\mathbf{x}) \right\}_{\lambda \in \Lambda}\right) \ge 0.
\]
\end{enumerate}
\end{corollary}

\bigskip

In the next result we extend the study for arbitrarily many quadratic inequalities of the finite quadratic program in \cite[Corollary 5.3]{jey-lee-li}, generalizing the Fritz John and Karush--Kuhn--Tucker theorems therein. Given a nonempty set $\Lambda$, we write $\ell^\infty (\Lambda)^*_+$ for the nonnegative cone of $\ell^\infty (\Lambda)_+$, that is, the set of all positive, linear and continuous functionals on $\ell^\infty (\Lambda)$. 

\bigskip

\begin{corollary}
Let $N \ge 1$, $X$ be a subset of $\mathbb{R}^N$ with $\mathbb{R}^N_+ \subset X$, $\mathbf{x}^0 \in X$ and $\Lambda $ be a nonempty set. Suppose that $\mathbf{A}  \in \mathbb{S}^N$,
$\mathbf{b}\in \mathbb{R}^N$,
$c\in \mathbb{R}$ satisfy that
\[
\left( \begin{array}{cc}
  \mathbf{A} & \mathbf{b} \\
  \mathbf{b}^T & 2 c
\end{array} \right)
\]
is a Z-matrix and that for all $\lambda \in \Lambda$,
$\mathbf{A}_\lambda  \in \mathbb{S}^N$,
$\mathbf{b}_\lambda \in \mathbb{R}^N$,
$c_\lambda \in \mathbb{R}$,
\[
\left( \begin{array}{cc}
  \mathbf{A}_\lambda & \mathbf{b}_\lambda \\
  \mathbf{b}_\lambda ^T & 2 c_\lambda
\end{array} \right)
\]
is also a Z-matrix. If, moreover,  $q, q_\lambda : \mathbb{R}^N \longrightarrow \mathbb{R}$ ($\lambda \in \Lambda$) are the quadratic functions given at each $\mathbf{x} \in \mathbb{R}^N$ by
\[
q(\mathbf{x}):=\frac{1}{2}\mathbf{x}^T \mathbf{A}\mathbf{x}+\mathbf{b}^T\mathbf{x}+c \quad \hbox{ and } \quad q_\lambda(\mathbf{x}):=\frac{1}{2}\mathbf{x}^T \mathbf{A}_\lambda\mathbf{x}+\mathbf{b}_\lambda^T\mathbf{x}+c_\lambda
\]
and satisfy that
\[
\mathbf{x} \in X \ \Rightarrow \ \left\{ q_\lambda (\mathbf{x}) \right\}_{\lambda \in \Lambda} \in \ell^\infty (\Lambda)
\]
and that the subset of $X$
\[
X_0:=\left\{ \mathbf{x} \in X: \ \sup_{\lambda \in \Lambda} q_\lambda (\mathbf{x}) \le 0 \right\}
\]
is nonempty, then we consider the quadratic program
\begin{equation*}\label{po:qp}\tag{QP}
\inf_{\mathbf{x} \in {X_0}} q(\mathbf{x}).
\end{equation*}
\begin{enumerate}
\item[{\rm (i)}] If $\mathbf{x}^0$ is an optimal solution of \eqref{po:qp}, then there exists $(y,L) \in \mathbb{R}_+ \times \ell^\infty (\Lambda)_+^*$ such that
\[
y+L(\mathbf{1})\neq 0
\]
and satisfying the Fritz John conditions:
\[
y q (\mathbf{x}) + L \left( \left\{ q_\lambda (\mathbf{x}) \right\}_{\lambda \in \Lambda} \right) \hbox{ attains its infimum on } X \hbox{ at } \mathbf{x}^0
\]
and
\[
L \left( \left\{ q_\lambda (\mathbf{x}^0) \right\}_{\lambda \in \Lambda} \right)=0.
\]
\item[{\rm (ii)}] Assume in addition that this Slater constraint qualification is fulfilled: there exists $\mathbf{x}^1 \in X$ with
\[
\sup_{\lambda \in \Lambda} q_\lambda (\mathbf{x}^1) <0.
\]
Then $\mathbf{x}^0$ is an optimal solution for \eqref{po:qp} if, and only if, there exists $L \in \ell^\infty (\Lambda)^*_+$ such that the following Karush--Kuhn--Tucker conditions are valid:
\[
q (\mathbf{x}) + L \left( \left\{ q_\lambda (\mathbf{x}) \right\}_{\lambda \in \Lambda} \right) \hbox{ attains its infimum on } X \hbox{ at } \mathbf{x}^0,
\]
\[
\sup_{\lambda \in \Lambda} q_\lambda (\mathbf{x}^0) \le 0,
\]
and
\[
L \left( \left\{ q_\lambda (\mathbf{x}^0) \right\}_{\lambda \in \Lambda} \right)=0.
\]

\end{enumerate}
\end{corollary}

\noindent \textsc{Proof.}
\fbox{(i)} Since $\mathbf{x}^0$ is an optimal solution of the quadratic program \eqref{po:qp}, then
\[
\inf_{\mathbf{x} \in X} \max \left\{ q(\mathbf{x})-q(\mathbf{x}^0), \sup_{\lambda \in \Lambda} q_\lambda (\mathbf{x}) \right\}=0,
\]
thus, the alternative in Corollary \ref{co:quadratic} provides a $(y,L) \in \mathbb{R}_+ \times \ell^\infty (\Lambda)^*_+$ with
\[
y+L(\mathbf{1})=1
\]
and such that
\[
\mathbf{x} \in X \ \Rightarrow \ 0 \le y(q(\mathbf{x})-q(\mathbf{x}^0))+ L\left(  \left\{ q_\lambda (\mathbf{x}) \right\}_{\lambda \in \Lambda}\right),
\]
i.e.,
\begin{equation}\label{eq:fj}
\mathbf{x} \in X \ \Rightarrow \ y q(\mathbf{x}^0) \le yq(\mathbf{x})+ L\left(  \left\{ q_\lambda (\mathbf{x}) \right\}_{\lambda \in \Lambda}\right).
\end{equation}
To conclude, it suffices to prove that
\[
0 = L\left(  \left\{ q_\lambda (\mathbf{x}^0) \right\}_{\lambda \in \Lambda}\right).
\]
But, on the one hand, inequality \eqref{eq:fj} for $\mathbf{x}^0$ implies
\[
0 \le L\left(  \left\{ q_\lambda (\mathbf{x}^0) \right\}_{\lambda \in \Lambda}\right);
\]
and, on the other hand, $\left\{ -q_\lambda(\mathbf{x}^0)\right\}_{\lambda \in \Lambda} \in \ell^\infty (\Lambda)_+$ ($\mathbf{x}^0$ being an optimal solution of \eqref{po:qp}, in particular, $\mathbf{x}^0 \in X_0$), which according to the positivity of $L$ yields
\[
 L\left(  \left\{ q_\lambda (\mathbf{x}^0) \right\}_{\lambda \in \Lambda}\right) \le 0.
\]

\noindent \fbox{(ii)} Let us first assume that $\mathbf{x}^0$ is an optimal solution for \eqref{po:qp}. Then, according to (i) we clearly have the announced statement (it suffices to divide by $y$, which is nonzero thanks to the Slater condition).

And conversely, given $\mathbf{x} \in X_0$, the assumptions clearly yield
\[
\begin{array}{rl}
    q(\mathbf{x}^0) & = \displaystyle q(\mathbf{x}^0) + L\left( \left\{ q(\mathbf{x}^0) \right\}_{\lambda \in \Lambda} \right) \\
                    & \le  q(\mathbf{x}) + L\left( \left\{ q(\mathbf{x}) \right\}_{\lambda \in \Lambda} \right) \\
                    & \le q(\mathbf{x})
\end{array}
\]
and since $\mathbf{x}^0 \in X_0$, then $\mathbf{x}^0$ is an optimal solution for the quadratic program \eqref{po:qp}.
\hfill$\textcolor[rgb]{0.00,0.00,0.60}{\Box}$

\section{Fenchel conjugate}\label{sect4}

\vspace{-0.5cm}

Before ending, we emphasize the applicability of the Gordan theorem with another result, now deriving a formula for the Fenchel conjugate of the supremum of a certain family of functions for which, unlike \cite[Proposition 5.3]{din-mo} or \cite[Remark 12.2]{bot}, we skip any topological assumption and weaken the convexity hypothesis.

If $E$ is a real topological vector space, $E^*$ denotes its
\textcolor[rgb]{0.89,0.39,0.02}{\textsl{topological dual space}}, i.e., the space of those linear and continuous functionals  $x^* : E
\longrightarrow \mathbb{R}$. Furthermore, given a function $f: E \longrightarrow \mathbb{R} \cup \{ +\infty
\}$, $f^*$ stands for its
\textcolor[rgb]{0.89,0.39,0.02}{\textsl{Fenchel conjugate}}, that
is, for each $x^* \in E^*$
\[
f^*(x^*) = \sup_{x \in E} (x^*(x)-f(x)).
\]

\bigskip

\begin{corollary}
If $E$ is a real topological vector space, $x_0^* \in E^*$, $\Lambda$ is a nonempty set and $\left\{ f_\lambda \right\}_{\lambda \in \Lambda}$ is a family of real-valued functions defined on $E$ such that
\[
\left\{ f_\lambda -x_0^* \right\}_{\lambda \in \Lambda} \hbox{ is infsup-convex on } E
\]
and
\[
x \in E \ \Rightarrow \ \left\{ f_\lambda (x) \right\}_{\lambda \in \Lambda} \in \ell^\infty (\Lambda),
\]
then
\[
\left( \sup_{\lambda \in \Lambda} f_\lambda \right)^*(x_0^*) = \min_{L \in \Delta_\Lambda} \left( L \left(  \left\{ f_\lambda ( \cdot ) \right\}_{\lambda \in \Lambda}
\right) \right)^*(x_0^*).
\]
In particular, if the functions $f_\lambda$ are convex, then the formula above holds for any $x_0^* \in E^*$.
\end{corollary}

\noindent \textsc{Proof.} Given $L \in \Delta_\Lambda $ and $x^* \in E^*$, we have that
\[
\begin{array}{rl}
     \displaystyle \left( \sup_{\lambda \in \Lambda} f_\lambda  \right)^*(x^*) & = \displaystyle  \sup_{x \in E} \left( x^*(x) -\sup_{\lambda \in \Lambda} f_\lambda f_\lambda (x) \right) \\
     & \le \displaystyle  \sup_{x \in E} \left( x^*(x) - L \left(  \left\{ f_\lambda (x) \right\}_{\lambda \in \Lambda}
\right) \right) \\
    & = \displaystyle \left( L \left(  \left\{ f_\lambda ( \cdot ) \right\}_{\lambda \in \Lambda}
\right) \right)^*(x^*),
\end{array}
\]
since
\[
\Phi \in \ell^\infty (\Lambda)  \  \Rightarrow \  L(\Phi) \le \sup_{\lambda \in \Lambda} \Phi (\lambda).
\]
Therefore, to finish the proof we must show that the other inequality is true for some $L \in \Delta_\Lambda$. So, let
\[
\alpha:= - \left( \sup_{\lambda \in \Lambda} f_\lambda \right)^*(x_0^*),
\]
which, by the previous inequality, can be assumed to be finite (otherwise, any $L \in \Delta_\Lambda$ satisfies the announced equality). Let us consider the family of functions
\[
h_\lambda (x):= -x_0^*(x)+f_\lambda (x) - \alpha , \qquad (x \in E, \ \lambda \in \Lambda).
\]
As
\[
\inf_{x \in E} \sup_{\lambda \in \Lambda} h_\lambda (x)=0,
\]
the Gordan theorem,  Theorem \ref{th:gordan}, and our hypothesis on the infsup-convexity of $\{f_\lambda -x_0^*\}_{\lambda \in \Lambda}$ on $E$, guarantee the existence of $L \in \Delta_\Lambda$ such that
\[
\inf_{x \in E} L \left( \left\{ h_\lambda (x) \right\}_{\lambda \in \Lambda} \right) \ge 0,
\]
that is ($L \in \Delta_\Lambda$),
\[
\inf_{x \in E} \left( -x_0^*(x) + L \left( \left\{ f_\lambda (x) \right\}_{\lambda \in \Lambda} \right) \right) \ge \alpha,
\]
or, in other words,
\[
\left( L \left(  \left\{ f_\lambda ( \cdot ) \right\}_{\lambda \in \Lambda} \right) \right)^*(x_0^*) \le \left( \sup_{\lambda \in \Lambda} f_\lambda \right)^*(x_0^*),
\]
as required.
\hfill$\textcolor[rgb]{0.00,0.00,0.60}{\Box}$

\bigskip

In particular, in view of Lemma \ref{le:quadratic} we conclude that if $N \ge 1$, $\Lambda$ is a nonemptyset, for all $\lambda \in \Lambda$,
$\mathbf{A}_\lambda  \in \mathbb{S}^N$,
$\mathbf{b}_\lambda \in \mathbb{R}^N$,
$c_\lambda \in \mathbb{R}$ in such a way that
\[
\left( \begin{array}{cc}
  \mathbf{A}_\lambda & \mathbf{b}_\lambda \\
  \mathbf{b}_\lambda ^T & 2 c_\lambda
\end{array} \right)
\]
is a Z-matrix, $q_\lambda : \mathbb{R}^N \longrightarrow \mathbb{R}$ is the quadratic function defined at each $\mathbf{x} \in \mathbb{R}^N$ by
\[
q_\lambda (\mathbf{x}):= \frac{1}{2}\mathbf{x}^T
\mathbf{A}_\lambda \mathbf{x} + \mathbf{b}_\lambda ^T \mathbf{x} +
c_\lambda
\]
with
\[
\mathbf{x} \in \mathbb{R}^N \ \Rightarrow \ \left\{ q_\lambda (\mathbf{x})\right\}_{\lambda \in \Lambda} \in \ell^\infty (\Lambda)
\]
and $\mathbf{x}_0 \in \mathbb{R}^N_+ $, then
\[
\left( \sup_{\lambda \in \Lambda} q_\lambda \right)^*(\mathbf{x}_0) = \min_{L \in \Delta_\Lambda} \left( L \left(  \left\{ q_\lambda ( \cdot ) \right\}_{\lambda \in \Lambda}
\right) \right)^*(\mathbf{x}_0).
\]

\section*{Acknowledgement}
\vspace{-0.5cm}

Research partially supported by project MTM2016-80676-P (AEI/FEDER, UE) and by Junta de Andaluc\'{\i}a Grant FQM359.

\end{document}